\documentclass[12pt]{article}%
\usepackage{amsfonts}
\usepackage{sw20bams}
\usepackage{amsmath}
\usepackage{amssymb}
\usepackage{graphicx}%
\providecommand{\U}[1]{\protect\rule{.1in}{.1in}}
\begin{document}

\title{Popa's \textquotedblleft Recurrent Sequences"\ and Reciprocity}
\author{Steven Finch}
\date{June 20, 2025}
\maketitle

\begin{abstract}
Dumitru Popa found asymptotic expansions for certain nonlinear recurrences,
but left open the numerical evaluation of associated constants. \ We address
this issue. \ A change of variables involving reciprocals and the algorithm of
Mavecha \&\ Laohakosol play a key role in our computations.

\end{abstract}

\footnotetext{Copyright \copyright \ 2025 by Steven R. Finch. All rights
reserved.}

If a positive real sequence $\{x_{k}\}_{k=0}^{\infty}$ satisfies%
\[%
\begin{array}
[c]{ccc}%
x_{k}\sim\alpha\,k+\beta\ln(k)+C, &  & \alpha\neq0
\end{array}
\]
as $k\rightarrow\infty$, then (by Proposition 12 of \cite{P1-popa} or
Proposition 7 of \cite{P2-popa})%
\[
y_{k}=\frac{1}{x_{k}}\,\sim\,\frac{1}{\alpha\,k}-\frac{\beta}{\alpha^{2}}%
\frac{\ln(k)}{k^{2}}-\frac{C}{\alpha^{2}}\frac{1}{k^{2}}\text{.}%
\]
Assuming $\{x_{k}\}$ satisfies the recurrence $x_{k+1}=f(x_{k})$, we have
$y_{k+1}=g(y_{k})$ where%
\[
g(y)=\frac{1}{f\left(  \dfrac{1}{y}\right)  }%
\]
because $1/y_{k+1}=f\left(  1/y_{k}\right)  $. \ Accurately computing the
constant $C$ requires many terms in the asymptotic expansion of $\{x_{k}\}$.
\ The work of de Bruijn \cite{dB-popa} and Bencherif \&\ Robin \cite{BR-popa}
gave rise to Mavecha \&\ Laohakosol's efficient algorithm \cite{ML-popa},
which cannot be applied directly to $\{x_{k}\}$ but (under some circumstances)
can be applied to $\{y_{k}\}$. \ Hence the interplay between $f(x)$~and $g(y)$
is crucial in our paper. \ We draw upon examples of Popa's \cite{P1-popa,
P2-popa}, as well as examples of our own \cite{F1-popa, F2-popa, F3-popa,
F4-popa, F5-popa}. \ We also make reference to Popa's addition theorem
(Theorems 3 / 4 of \cite{P1-popa} / \cite{P2-popa}) involving%
\[%
\begin{array}
[c]{ccc}%
f(x)=x+\varphi\left(  \dfrac{1}{x}\right)  &  & \text{for smooth }%
\varphi(x)>0\,\,\forall x\geq0
\end{array}
\]
and Popa's multiplication theorem (Theorems 6 / 5 of \cite{P1-popa} /
\cite{P2-popa}) involving%
\[%
\begin{array}
[c]{ccc}%
f(x)=x\cdot\psi\left(  \dfrac{1}{x}\right)  &  & \text{for smooth }%
\psi(x)>1\,\,\forall x>0\text{; }\psi(0)=1\text{; }\psi^{\prime}%
(0)\neq0\text{.}%
\end{array}
\]
For simplicity, initial conditions $x_{0}=y_{0}=1$ are presumed throughout.

\section{$p$-Sequences}

Let $f(x)=x^{1-p}(1+x)^{p}$ where $0<p<1$, then $g(y)=y/(1+y)^{p}$. \ By the
multiplication theorem with $\psi(x)=(1+x)^{p}$,%
\[
x_{k}\sim\alpha\,k+\beta\ln(k)+C+\gamma\frac{\ln(k)}{k}+\delta\,\frac{1}{k}%
\]
where%
\[%
\begin{array}
[c]{ccccccc}%
\alpha=p, &  & \beta=-\dfrac{1}{2}+\dfrac{p}{2}, &  & \gamma=\dfrac{1}%
{4p}-\dfrac{1}{2}+\dfrac{p}{4}, &  & \delta=-\dfrac{1}{12p}-\dfrac{C}%
{2p}+\dfrac{1}{4}+\dfrac{C}{2}-\dfrac{p}{6}.
\end{array}
\]
We cover the cases $p=1/2,1/3$ and $2/3$ separately. \ The boundary case $p=1$
is trivial. \ There is otherwise no reason to restrict $p<1$; we deal with
$p=2$ as well.

\subsection{Case $p=1/2$}

This occurred as Corollary 9 in \cite{P1-popa} and in Section 5 of
\cite{F3-popa} (for the latter, a brute-force method was employed to calculate
$C$ rather than the Mavecha-Laohakosol algorithm). \ We find
\[
x_{k}\sim\frac{1}{2}k-\frac{1}{4}\ln(k)+C+\frac{1}{8}\frac{\ln(k)}{k}-\frac
{C}{2}\frac{1}{k},
\]%
\begin{align*}
y_{k}  &  \sim\frac{2}{k}+\frac{\ln(k)}{k^{2}}-\frac{4C}{k^{2}}+\frac{1}%
{2}\frac{\ln(k)^{2}}{k^{3}}-\left(  \frac{1}{2}+4C\right)  \frac{\ln(k)}%
{k^{3}}+\left(  2C+8C^{2}\right)  \frac{1}{k^{3}}+\frac{1}{4}\frac{\ln(k)^{3}%
}{k^{4}}\\
&  -\left(  \frac{5}{8}+3C\right)  \frac{\ln(k)^{2}}{k^{4}}+\left(  \frac
{1}{4}+5C+12C^{2}\right)  \frac{\ln(k)}{k^{4}}-\left(  -\frac{1}{24}%
+C+10C^{2}+16C^{3}\right)  \frac{1}{k^{4}}%
\end{align*}
and thus%
\[
-\lim_{k\rightarrow\infty}k^{2}\left(  y_{k}-\frac{2}{k}-\frac{\ln(k)}{k^{2}%
}\right)  =4C,
\]%
\[
\lim_{k\rightarrow\infty}\left(  x_{k}-\frac{1}{2}k+\frac{1}{4}\ln(k)\right)
=C=1.1751774424585571398132856....
\]

\subsection{Case $p=1/3$}

We find%
\[
x_{k}\sim\frac{1}{3}k-\frac{1}{3}\ln(k)+C+\frac{1}{3}\frac{\ln(k)}{k}-\left(
\frac{1}{18}+C\right)  \frac{1}{k},
\]%
\begin{align*}
y_{k}  &  \sim\frac{3}{k}+3\frac{\ln(k)}{k^{2}}-\frac{9C}{k^{2}}+3\frac
{\ln(k)^{2}}{k^{3}}-\left(  3+18C\right)  \frac{\ln(k)}{k^{3}}+\left(
\frac{1}{2}+9C+27C^{2}\right)  \frac{1}{k^{3}}+3\frac{\ln(k)^{3}}{k^{4}}\\
&  -\left(  \frac{15}{2}+27C\right)  \frac{\ln(k)^{2}}{k^{4}}+\left(  \frac
{9}{2}+45C+81C^{2}\right)  \frac{\ln(k)}{k^{4}}-\left(  \frac{1}{2}+\frac
{27C}{2}+\frac{135C^{2}}{2}+81C^{3}\right)  \frac{1}{k^{4}}%
\end{align*}
and thus
\[
-\lim_{k\rightarrow\infty}k^{2}\left(  y_{k}-\frac{3}{k}-3\frac{\ln(k)}{k^{2}%
}\right)  =9C,
\]%
\[
\lim_{k\rightarrow\infty}\left(  x_{k}-\frac{1}{3}k+\frac{1}{3}\ln(k)\right)
=C=1.3842423952717873718895461....
\]
\qquad

\subsection{Case $p=2/3$}

We find%
\[
x_{k}\sim\frac{2}{3}k-\frac{1}{6}\ln(k)+C+\frac{1}{24}\frac{\ln(k)}{k}-\left(
-\frac{1}{72}+\frac{C}{4}\right)  \frac{1}{k},
\]%
\begin{align*}
y_{k}  &  \sim\frac{3}{2k}+\frac{3}{8}\frac{\ln(k)}{k^{2}}-\frac{9C}{4}%
\frac{1}{k^{2}}+\frac{3}{32}\frac{\ln(k)^{2}}{k^{3}}-\left(  \frac{3}%
{32}+\frac{9C}{8}\right)  \frac{\ln(k)}{k^{3}}+\left(  -\frac{1}{32}+\frac
{9C}{16}+\frac{27C^{2}}{8}\right)  \frac{1}{k^{3}}\\
&  +\frac{3}{128}\frac{\ln(k)^{3}}{k^{4}}-\left(  \frac{15}{256}+\frac
{27C}{64}\right)  \frac{\ln(k)^{2}}{k^{4}}+\left(  \frac{45C}{64}%
+\frac{81C^{2}}{32}\right)  \frac{\ln(k)}{k^{4}}\\
&  -\left(  -\frac{5}{256}+\frac{135C^{2}}{64}+\frac{81C^{3}}{16}\right)
\frac{1}{k^{4}}%
\end{align*}
and thus
\[
-\lim_{k\rightarrow\infty}k^{2}\left(  y_{k}-\frac{3}{2k}-\frac{3}{8}\frac
{\ln(k)}{k^{2}}\right)  =\frac{9C}{4},
\]%
\[
\lim_{k\rightarrow\infty}\left(  x_{k}-\frac{2}{3}k+\frac{1}{6}\ln(k)\right)
=C=1.0603463553715904094868689....
\]

\subsection{Case $p=2$}

We find%
\[
x_{k}\sim2k+\frac{1}{2}\ln(k)+C+\frac{1}{8}\frac{\ln(k)}{k}+\left(  -\frac
{1}{8}+\frac{C}{4}\right)  \frac{1}{k},
\]%
\begin{align*}
y_{k}  &  \sim\frac{1}{2k}-\frac{1}{8}\frac{\ln(k)}{k^{2}}-\frac{C}{4}\frac
{1}{k^{2}}+\frac{1}{32}\frac{\ln(k)^{2}}{k^{3}}+\left(  -\frac{1}{32}+\frac
{C}{8}\right)  \frac{\ln(k)}{k^{3}}+\left(  \frac{1}{32}-\frac{C}{16}%
+\frac{C^{2}}{8}\right)  \frac{1}{k^{3}}\\
&  -\frac{1}{128}\frac{\ln(k)^{3}}{k^{4}}+\left(  \frac{5}{256}-\frac{3C}%
{64}\right)  \frac{\ln(k)^{2}}{k^{4}}-\left(  \frac{1}{32}-\frac{5C}{64}%
+\frac{3C^{2}}{32}\right)  \frac{\ln(k)}{k^{4}}\\
&  +\left(  \frac{11}{768}-\frac{C}{16}+\frac{5C^{2}}{64}-\frac{C^{3}}%
{16}\right)  \frac{1}{k^{4}}%
\end{align*}
and thus
\[
-\lim_{k\rightarrow\infty}k^{2}\left(  y_{k}-\frac{1}{2k}+\frac{1}{8}\frac
{\ln(k)}{k^{2}}\right)  =\frac{C}{4},
\]%
\[
\lim_{k\rightarrow\infty}\left(  x_{k}-2k-\frac{1}{2}\ln(k)\right)
=C=1.7231423751374234610635627....
\]

\section{Radicals}

\subsection{Case $\sqrt{1+x+x^{2}}$}

Corollary 9 of \cite{P1-popa} was devoted to $\sqrt{x^{2}+a\,x+b}$ where $a>0$
and $b\geq0$. \ We examined already $a=1$ and $b=0$ in Section 1.1. \ The
boundary case $a=0$ and $b=1$ is trivial. \ A vital participant here is the
multiplication theorem with $\psi(x)=\sqrt{1+a\,x+b\,x^{2}}$. \ For $a=b=1$,
we find%
\[%
\begin{array}
[c]{ccc}%
f(x)=\sqrt{1+x+x^{2}}, &  & g(y)=\dfrac{y}{\sqrt{1+y+y^{2}}},
\end{array}
\]%
\[
x_{k}\sim\frac{1}{2}k+\frac{3}{4}\ln(k)+C+\frac{9}{8}\frac{\ln(k)}{k}+\left(
\frac{3}{2}+\frac{3C}{2}\right)  \frac{1}{k},
\]%
\begin{align*}
y_{k}  &  \sim\frac{2}{k}-3\frac{\ln(k)}{k^{2}}-\frac{4C}{k^{2}}+\frac{9}%
{2}\frac{\ln(k)^{2}}{k^{3}}+\left(  -\frac{9}{2}+12C\right)  \frac{\ln
(k)}{k^{3}}-\left(  6+6C-8C^{2}\right)  \frac{1}{k^{3}}\\
&  -\frac{27}{4}\frac{\ln(k)^{3}}{k^{4}}+\left(  \frac{135}{8}-27C\right)
\frac{\ln(k)^{2}}{k^{4}}+\left(  \frac{81}{4}+45C-36C^{2}\right)  \frac
{\ln(k)}{k^{4}}\\
&  +\left(  -\frac{25}{8}+27C+30C^{2}-16C^{3}\right)  \frac{1}{k^{4}}%
\end{align*}
and thus
\[
-\lim_{k\rightarrow\infty}k^{2}\left(  y_{k}-\frac{2}{k}+3\frac{\ln(k)}{k^{2}%
}\right)  =4C,
\]%
\[
\lim_{k\rightarrow\infty}\left(  x_{k}-\frac{1}{2}k-\frac{3}{4}\ln(k)\right)
=C=-0.0431693967745555915829486....
\]
Incidentally, Proposition 13 of \cite{P1-popa} was devoted to $y/\sqrt
{1+a\,y+b\,y^{2}}$; Popa was indeed aware of reciprocity. \ He did not,
however, pursue its use in calculating $C$. \ Propositions 10 \&\ 14
demonstrated likewise, but we omit these for reasons of space.

\subsection{Case $x+\sqrt{1+\frac{1}{x}}$}

Corollary 5 of \cite{P1-popa} was devoted to $x+\sqrt{a+\frac{b}{x}}$ where
$a>0$ and $b>0$. \ A vital participant here is the addition theorem with
$\varphi(x)=\sqrt{a+b\,x}$. \ For $a=b=1$, we find%
\[%
\begin{array}
[c]{ccc}%
f(x)=x+\sqrt{1+\dfrac{1}{x}}, &  & g(y)=\dfrac{y}{1+y\sqrt{1+y}},
\end{array}
\]%
\[
x_{k}\sim k+\frac{1}{2}\ln(k)+C+\frac{1}{4}\frac{\ln(k)}{k}+\left(  \frac
{1}{8}+\frac{C}{2}\right)  \frac{1}{k},
\]%
\begin{align*}
y_{k}  &  \sim\frac{1}{k}-\frac{1}{2}\frac{\ln(k)}{k^{2}}-\frac{C}{k^{2}%
}+\frac{1}{4}\frac{\ln(k)^{2}}{k^{3}}+\left(  -\frac{1}{4}+C\right)  \frac
{\ln(k)}{k^{3}}-\left(  \frac{1}{8}+\frac{C}{2}-C^{2}\right)  \frac{1}{k^{3}%
}\\
&  -\frac{1}{8}\frac{\ln(k)^{3}}{k^{4}}+\left(  \frac{5}{16}-\frac{3C}%
{4}\right)  \frac{\ln(k)^{2}}{k^{4}}+\left(  \frac{1}{16}+\frac{5C}{4}%
-\frac{3C^{2}}{2}\right)  \frac{\ln(k)}{k^{4}}\\
&  +\left(  \frac{1}{96}+\frac{C}{8}+\frac{5C^{2}}{4}-C^{3}\right)  \frac
{1}{k^{4}}%
\end{align*}
and thus
\[
-\lim_{k\rightarrow\infty}k^{2}\left(  y_{k}-\frac{1}{k}+\frac{1}{2}\frac
{\ln(k)}{k^{2}}\right)  =C,
\]%
\[
\lim_{k\rightarrow\infty}\left(  x_{k}-k-\frac{1}{2}\ln(k)\right)
=C=0.9330050241078502218961438....
\]

\section{Exponentials / Logarithms}

It is surprising that the first \& second examples below do not appear in
\cite{P1-popa}, due to their simplicity. \ The seventh \&\ eighth examples are
understandably missing, as they are more complicated.

\subsection{Case $x\exp\left(  \frac{1}{x}\right)  $}

We find, via multiplication,%
\[%
\begin{array}
[c]{ccc}%
f(x)=x\exp\left(  \dfrac{1}{x}\right)  , &  & g(y)=y\exp(-y),
\end{array}
\]%
\[
x_{k}\sim k+\frac{1}{2}\ln(k)+C+\frac{1}{4}\frac{\ln(k)}{k}-\left(  \frac
{1}{6}-\frac{C}{2}\right)  \frac{1}{k},
\]%
\begin{align*}
y_{k}  &  \sim\frac{1}{k}-\frac{1}{2}\frac{\ln(k)}{k^{2}}-\frac{C}{k^{2}%
}+\frac{1}{4}\frac{\ln(k)^{2}}{k^{3}}+\left(  -\frac{1}{4}+C\right)  \frac
{\ln(k)}{k^{3}}\\
&  +\left(  \frac{1}{6}-\frac{C}{2}+C^{2}\right)  \frac{1}{k^{3}}-\frac{1}%
{8}\frac{\ln(k)^{3}}{k^{4}}+\left(  \frac{5}{16}-\frac{3C}{4}\right)
\frac{\ln(k)^{2}}{k^{4}}\\
&  -\left(  \frac{3}{8}-\frac{5C}{4}+\frac{3C^{2}}{2}\right)  \frac{\ln
(k)}{k^{4}}+\left(  \frac{7}{48}-\frac{3C}{4}+\frac{5C^{2}}{4}-C^{3}\right)
\frac{1}{k^{4}}%
\end{align*}
and thus%
\[
-\lim_{k\rightarrow\infty}k^{2}\left(  y_{k}-\frac{1}{k}+\frac{1}{2}\frac
{\ln(k)}{k^{2}}\right)  =C,
\]%
\[
\lim_{k\rightarrow\infty}\left(  x_{k}-k-\frac{1}{2}\ln(k)\right)
=C=1.2902472086877642916676156....
\]
The asymptotic series for $\{y_{k}\}$ appeared in Section 14 of \cite{F4-popa}%
; the map $y\mapsto-g(-y)$ is the functional inverse of the Lambert map
$y\mapsto W(y)$ and hence $g$ \&\ $W$ are kindred functions (i.e.,
corresponding recursions possess remarkably similar expansions).

\subsection{Case $x+\exp\left(  \frac{1}{x}\right)  $}

We find, via addition,%
\[%
\begin{array}
[c]{ccc}%
f(x)=x+\exp\left(  \dfrac{1}{x}\right)  , &  & g(y)=\dfrac{y}{1+y\exp(y)},
\end{array}
\]%
\[
x_{k}\sim k+\ln(k)+C+\frac{\ln(k)}{k}+\frac{C}{k},
\]%
\begin{align*}
y_{k}  &  \sim\frac{1}{k}-\frac{\ln(k)}{k^{2}}-\frac{C}{k^{2}}+\frac
{\ln(k)^{2}}{k^{3}}+\left(  -1+2C\right)  \frac{\ln(k)}{k^{3}}+\left(
-C+C^{2}\right)  \frac{1}{k^{3}}\\
&  -\frac{\ln(k)^{3}}{k^{4}}+\left(  \frac{5}{2}-3C\right)  \frac{\ln(k)^{2}%
}{k^{4}}+\left(  -1+5C-3C^{2}\right)  \frac{\ln(k)}{k^{4}}\\
&  +\left(  \frac{1}{6}-C+\frac{5C^{2}}{2}-C^{3}\right)  \frac{1}{k^{4}}%
\end{align*}
and thus%
\[
-\lim_{k\rightarrow\infty}k^{2}\left(  y_{k}-\frac{1}{k}+\frac{\ln(k)}{k^{2}%
}\right)  =C,
\]%
\[
\lim_{k\rightarrow\infty}\left(  x_{k}-k-\ln(k)\right)
=C=1.4196070070201119940555181....
\]

\subsection{Case $x\left[  1+\ln\left(  1+\frac{1}{x}\right)  \right]  $}

This occurred as Corollary 8 in \cite{P1-popa}, to be discussed more later.
\ We find, via multiplication,%
\[%
\begin{array}
[c]{ccc}%
f(x)=x\left[  1+\ln\left(  1+\dfrac{1}{x}\right)  \right]  , &  &
g(y)=\dfrac{y}{1+\ln(1+y)},
\end{array}
\]%
\[
x_{k}\sim k-\frac{1}{2}\ln(k)+C+\frac{1}{4}\frac{\ln(k)}{k}+\left(  \frac
{1}{6}-\frac{C}{2}\right)  \frac{1}{k},
\]%
\begin{align*}
y_{k}  &  \sim\frac{1}{k}+\frac{1}{2}\frac{\ln(k)}{k^{2}}-\frac{C}{k^{2}%
}+\frac{1}{4}\frac{\ln(k)^{2}}{k^{3}}-\left(  \frac{1}{4}+C\right)  \frac
{\ln(k)}{k^{3}}\\
&  +\left(  -\frac{1}{6}+\frac{C}{2}+C^{2}\right)  \frac{1}{k^{3}}+\frac{1}%
{8}\frac{\ln(k)^{3}}{k^{4}}-\left(  \frac{5}{16}+\frac{3C}{4}\right)
\frac{\ln(k)^{2}}{k^{4}}\\
&  +\left(  -\frac{1}{8}+\frac{5C}{4}+\frac{3C^{2}}{2}\right)  \frac{\ln
(k)}{k^{4}}+\left(  \frac{1}{8}+\frac{C}{4}-\frac{5C^{2}}{4}-C^{3}\right)
\frac{1}{k^{4}}%
\end{align*}
and thus%
\[
-\lim_{k\rightarrow\infty}k^{2}\left(  y_{k}-\frac{1}{k}-\frac{1}{2}\frac
{\ln(k)}{k^{2}}\right)  =C,
\]%
\[
\lim_{k\rightarrow\infty}\left(  x_{k}-k+\frac{1}{2}\ln(k)\right)
=C=0.8725870124185800733516473....
\]

\subsection{Case $x+1+\ln\left(  1+\frac{1}{x}\right)  $}

Corollary 4 of \cite{P1-popa} was devoted to $x+\ln\left(  a+\frac{b}%
{x}\right)  $ where $a>1$ and $b>0$. We set $a=b=e$ and find, via addition,%
\[%
\begin{array}
[c]{ccc}%
f(x)=x+1+\ln\left(  1+\dfrac{1}{x}\right)  , &  & g(y)=\dfrac{y}%
{y+1+y\ln(1+y)},
\end{array}
\]%
\[
x_{k}\sim k+\ln(k)+C+\frac{\ln(k)}{k}+\left(  1+C\right)  \frac{1}{k},
\]%
\begin{align*}
y_{k}  &  \sim\frac{1}{k}-\frac{\ln(k)}{k^{2}}-\frac{C}{k^{2}}+\frac
{\ln(k)^{2}}{k^{3}}+\left(  -1+2C\right)  \frac{\ln(k)}{k^{3}}+\left(
-1-C+C^{2}\right)  \frac{1}{k^{3}}\\
&  -\frac{\ln(k)^{3}}{k^{4}}+\left(  \frac{5}{2}-3C\right)  \frac{\ln(k)^{2}%
}{k^{4}}+\left(  2+5C-3C^{2}\right)  \frac{\ln(k)}{k^{4}}\\
&  +\left(  -\frac{1}{4}+2C+\frac{5C^{2}}{2}-C^{3}\right)  \frac{1}{k^{4}}%
\end{align*}
and thus%
\[
-\lim_{k\rightarrow\infty}k^{2}\left(  y_{k}-\frac{1}{k}+\frac{\ln(k)}{k^{2}%
}\right)  =C,
\]%
\[
\lim_{k\rightarrow\infty}\left(  x_{k}-k-\ln(k)\right)
=C=0.409139256166675934818465....
\]

\subsection{Case $x\exp\left(  \sqrt{1+\frac{1}{x}}-1\right)  $}

This occurred as Corollary 7 in \cite{P1-popa}, to be discussed more later.
\ We find, via multiplication,%
\[%
\begin{array}
[c]{ccc}%
f(x)=x\exp\left(  \sqrt{1+\dfrac{1}{x}}-1\right)  , &  & g(y)=y\exp\left(
1-\sqrt{1+y}\right)  ,
\end{array}
\]%
\[
x_{k}\sim\frac{k}{2}+C-\frac{1}{12k},
\]%
\begin{align*}
y_{k}  &  \sim\frac{2}{k}-\frac{4C}{k^{2}}+\left(  \frac{1}{3}+8C^{2}\right)
\frac{1}{k^{3}}-\left(  \frac{1}{24}+2C+16C^{3}\right)  \frac{1}{k^{4}%
}+\left(  \frac{151}{1080}+\frac{C}{3}+8C^{2}+32C^{4}\right)  \frac{1}{k^{5}%
}\\
&  -\left(  \frac{67}{1440}+\frac{151C}{108}+\frac{5C^{2}}{3}+\frac{80C^{3}%
}{3}+64C^{5}\right)  \frac{1}{k^{6}}%
\end{align*}
and thus%
\[
-\lim_{k\rightarrow\infty}k^{2}\left(  y_{k}-\frac{2}{k}\right)  =4C,
\]%
\[
\lim_{k\rightarrow\infty}\left(  x_{k}-\frac{k}{2}\right)
=C=1.0398066960295413916327213....
\]
These are the first asymptotic expansions in this paper for which all
logarithmic terms evidently vanish.

\subsection{Case $x\left[  1+\ln\left(  1+\frac{2}{x}+\frac{2}{x^{2}}\right)
\right]  $}

This occurred as Corollary 8 in \cite{P1-popa}, to be discussed more later.
\ Via multiplication,%
\[%
\begin{array}
[c]{ccc}%
f(x)=x\left[  1+\ln\left(  1+\dfrac{2}{x}+\dfrac{2}{x^{2}}\right)  \right]
, &  & g(y)=\dfrac{y}{1+\ln(1+2y+2y^{2})},
\end{array}
\]%
\[
x_{k}\sim2k+C+\frac{1}{3k},
\]%
\begin{align*}
y_{k}  &  \sim\frac{1}{2k}-\frac{C}{4k^{2}}+\left(  -\frac{1}{12}+\frac{C^{2}%
}{8}\right)  \frac{1}{k^{3}}-\left(  \frac{1}{96}-\frac{C}{8}+\frac{C^{3}}%
{16}\right)  \frac{1}{k^{4}}+\left(  \frac{139}{4320}+\frac{C}{48}-\frac
{C^{2}}{8}+\frac{C^{4}}{32}\right)  \frac{1}{k^{5}}\\
&  +\left(  \frac{13}{960}-\frac{139C}{1728}-\frac{5C^{2}}{192}+\frac{5C^{3}%
}{48}-\frac{C^{5}}{64}\right)  \frac{1}{k^{6}}%
\end{align*}
and thus%
\[
-\lim_{k\rightarrow\infty}k^{2}\left(  y_{k}-\frac{1}{2k}\right)  =\frac{C}%
{4},
\]%
\[
\lim_{k\rightarrow\infty}\left(  x_{k}-2k\right)
=C=0.3410952769697805797926024....
\]
Like the asymptotic expansions in Section 3.5, all logarithmic terms evidently
vanish and convergence is quick.

\subsection{Case $x\exp\left(  \sqrt{\frac{1}{x}}\right)  $}

Corollary 7 of \cite{P1-popa} was devoted to $x\exp\left(  \sqrt{a^{2}%
+\frac{1}{x}}-a\right)  $ where $a>0$; the case $a=1$ was discussed earlier in
Section 3.5. \ The case $a=0$ is on the boundary of allowable values and
Popa's expansion for $x_{k}$ does not apply to $f(x)=x\exp\left(  \sqrt
{\frac{1}{x}}\right)  $. \ Mavecha \&\ Laohakosol's algorithm \textit{does},
however, apply to $g(y)=y\exp\left(  -\sqrt{y}\right)  $:%
\begin{align*}
y_{k}  &  \sim\frac{4}{k^{2}}-4\frac{\ln(k)}{k^{3}}-\frac{16C}{k^{3}}%
+3\frac{\ln(k)^{2}}{k^{4}}+\left(  -2+24C\right)  \frac{\ln(k)}{k^{4}}+\left(
\frac{4}{3}-8C+48C^{2}\right)  \frac{1}{k^{4}}\\
&  -2\frac{\ln(k)^{3}}{k^{5}}+\left(  \frac{7}{2}-24C\right)  \frac{\ln
(k)^{2}}{k^{5}}+\left(  -\frac{11}{3}+28C-96C^{2}\right)  \frac{\ln(k)}{k^{5}%
}\\
&  +\left(  \frac{7}{6}-\frac{44C}{3}+56C^{2}-128C^{3}\right)  \frac{1}{k^{5}}%
\end{align*}
and thus%
\[
-\lim_{k\rightarrow\infty}k^{3}\left(  y_{k}-\frac{4}{k^{2}}+4\frac{\ln
(k)}{k^{3}}\right)  =16C,
\]%
\[
\lim_{k\rightarrow\infty}\frac{1}{k}\left(  x_{k}-\frac{1}{4}k^{2}-\frac{1}%
{4}k\ln(k)\right)  =C=0.8791712792948618603132189....
\]

\subsection{Case $x\left[  1+\ln\left(  1+\frac{1}{x^{2}}\right)  \right]  $}

Corollary 8 of \cite{P1-popa} was devoted to $x\ln\left(  e+\frac{a}{x}%
+\frac{b}{x^{2}}\right)  $ where $a>0$ and $b\geq0$; the case $a=e$ \&\ $b=0$
was discussed earlier in Section 3.3 and the case $a=b=2e$ in Section
3.6.\ \ The case $a=0$ \&\ $b=e$ is on the boundary of allowable values and
Popa's expansion for $x_{k}$ does not apply to $f(x)=x\left[  1+\ln\left(
1+\frac{1}{x^{2}}\right)  \right]  $. \ Mavecha \&\ Laohakosol's algorithm
\textit{does}, however, apply to $g(y)=\frac{y}{1+\ln\left(  1+y^{2}\right)
}$:%
\begin{align*}
\sqrt{2}y_{k}  &  \sim\frac{1}{k^{1/2}}-\frac{C}{k^{3/2}}+\left(  -\frac
{1}{48}+\frac{3C^{2}}{2}\right)  \frac{1}{k^{5/2}}+\left(  -\frac{1}%
{256}+\frac{5C}{48}-\frac{5C^{3}}{2}\right)  \frac{1}{k^{7/2}}\\
&  +\left(  \frac{139}{69120}+\frac{7C}{256}-\frac{35C^{2}}{96}+\frac{35C^{4}%
}{8}\right)  \frac{1}{k^{9/2}}\\
&  +\left(  \frac{55}{36864}-\frac{139C}{7680}-\frac{63C^{2}}{512}%
+\frac{35C^{3}}{32}-\frac{63C^{5}}{8}\right)  \frac{1}{k^{11/2}}%
\end{align*}
and thus%
\[
-\lim_{k\rightarrow\infty}k^{3/2}\left(  \sqrt{2}y_{k}-\frac{1}{k^{1/2}%
}\right)  =C,
\]%
\[
\lim_{k\rightarrow\infty}k^{1/2}\left(  \frac{1}{\sqrt{2}}x_{k}-k^{1/2}%
\right)  =C=0.2005534003696638830775944....
\]

\section{$q$-Sequences}

Let $f(x)=x+1/x^{q-1}$ where $q>1$, then $g(y)=y/(1+y^{q})$. \ Details of a
transformation, followed by the multiplication theorem with $\psi
(x)=(1+x)^{q}$, were given in Section 3 of \cite{F5-popa}. \ This yielded a
six-term series, the first three terms of which are
\[
q^{1-1/q}\,x_{k}\sim q\,k^{1/q}+\left(  -\frac{1}{2q}+\frac{1}{2}\right)
\frac{\ln(k)}{k^{1-1/q}}+\frac{C}{k^{1-1/q}}.
\]
We cover the cases $q=2,3,3/2,4$ and $4/3$ separately. \ The boundary case
$q=1$ is trivial. \ There is otherwise no reason to restrict $q>1$; we deal
with $q=1/2$ as well.

\subsection{Case $q=2$}

The asymptotics of $x_{k+1}=x_{k}+1/x_{k}$ have generated considerable
interest \cite{Nw-popa, BL-popa, Mnt-popa, LPD-popa, Slv-popa, KAL-popa};
Corollaries 15 / 8 of \cite{P1-popa} / \cite{P2-popa} and Sections 4 / 3 of
\cite{F2-popa} / \cite{F5-popa} also. \ We find%
\begin{align*}
2^{1/2}y_{k}  &  \sim\frac{1}{k^{1/2}}-\frac{1}{8}\frac{\ln(k)}{k^{3/2}}%
-\frac{C}{2}\frac{1}{k^{3/2}}+\frac{3}{128}\frac{\ln(k)^{2}}{k^{5/2}}+\left(
-\frac{1}{32}+\frac{3C}{16}\right)  \frac{\ln(k)}{k^{5/2}}\\
&  +\left(  \frac{1}{32}-\frac{C}{8}+\frac{3C^{2}}{8}\right)  \frac{1}%
{k^{5/2}}-\frac{5}{1024}\frac{\ln(k)^{3}}{k^{7/2}}+\left(  \frac{1}{64}%
-\frac{15C}{256}\right)  \frac{\ln(k)^{2}}{k^{7/2}}\\
&  +\left(  -\frac{7}{256}+\frac{C}{8}-\frac{15C^{2}}{64}\right)  \frac
{\ln(k)}{k^{7/2}}+\left(  \frac{11}{768}-\frac{7C}{64}+\frac{C^{2}}{4}%
-\frac{5C^{3}}{16}\right)  \frac{1}{k^{7/2}}%
\end{align*}
and thus
\[
-\lim_{k\rightarrow\infty}k^{3/2}\left(  2^{1/2}y_{k}-\frac{1}{k^{1/2}}%
+\frac{1}{8}\frac{\ln(k)}{k^{3/2}}\right)  =\frac{C}{2},
\]%
\[
\lim_{k\rightarrow\infty}k^{1/2}\left(  2^{1/2}x_{k}-2k^{1/2}-\frac{1}{4}%
\frac{\ln(k)}{k^{1/2}}\right)  =C=0.8615711875687117305317813....
\]
An apt\ name \textquotedblleft add-the-reciprocal sequence\textquotedblright%
\ has been proposed for $1,2,5/2,29/10,941/290,...$ \cite{SS-popa, FS-popa}.

\subsection{Case $q=3$}

We find%
\begin{align*}
3^{1/3}y_{k}  &  \sim\frac{1}{k^{1/3}}-\frac{1}{9}\frac{\ln(k)}{k^{4/3}}%
-\frac{C}{3}\frac{1}{k^{4/3}}+\frac{2}{81}\frac{\ln(k)^{2}}{k^{7/3}}+\left(
-\frac{1}{27}+\frac{4C}{27}\right)  \frac{\ln(k)}{k^{7/3}}\\
&  +\left(  \frac{5}{162}-\frac{C}{9}+\frac{2C^{2}}{9}\right)  \frac
{1}{k^{7/3}}-\frac{14}{2187}\frac{\ln(k)^{3}}{k^{10/3}}+\left(  \frac{11}%
{486}-\frac{14C}{243}\right)  \frac{\ln(k)^{2}}{k^{10/3}}\\
&  +\left(  -\frac{53}{1458}+\frac{11C}{81}-\frac{14C^{2}}{81}\right)
\frac{\ln(k)}{k^{10/3}}+\left(  \frac{1}{54}-\frac{53C}{486}+\frac{11C^{2}%
}{54}-\frac{14C^{3}}{81}\right)  \frac{1}{k^{10/3}}%
\end{align*}
and thus
\[
-\lim_{k\rightarrow\infty}k^{4/3}\left(  3^{1/3}y_{k}-\frac{1}{k^{1/3}}%
+\frac{1}{9}\frac{\ln(k)}{k^{4/3}}\right)  =\frac{C}{3},
\]%
\[
\lim_{k\rightarrow\infty}k^{2/3}\left(  3^{2/3}x_{k}-3k^{1/3}-\frac{1}{3}%
\frac{\ln(k)}{k^{2/3}}\right)  =C=1.3784186157718345713984647....
\]

\subsection{Case $q=3/2$}

We find%
\begin{align*}
\left(  \frac{3}{2}\right)  ^{2/3}y_{k}  &  \sim\frac{1}{k^{2/3}}-\frac{1}%
{9}\frac{\ln(k)}{k^{5/3}}-\frac{2C}{3}\frac{1}{k^{5/3}}+\frac{5}{324}\frac
{\ln(k)^{2}}{k^{8/3}}+\left(  -\frac{1}{54}+\frac{5C}{27}\right)  \frac
{\ln(k)}{k^{8/3}}\\
&  +\left(  \frac{2}{81}-\frac{C}{9}+\frac{5C^{2}}{9}\right)  \frac{1}%
{k^{8/3}}-\frac{5}{2187}\frac{\ln(k)^{3}}{k^{11/3}}+\left(  \frac{13}%
{1944}-\frac{10C}{243}\right)  \frac{\ln(k)^{2}}{k^{11/3}}\\
&  +\left(  -\frac{41}{2916}+\frac{13C}{162}-\frac{20C^{2}}{81}\right)
\frac{\ln(k)}{k^{11/3}}+\left(  \frac{5}{648}-\frac{41C}{486}+\frac{13C^{2}%
}{54}-\frac{40C^{3}}{81}\right)  \frac{1}{k^{11/3}}%
\end{align*}
and thus
\[
-\lim_{k\rightarrow\infty}k^{5/3}\left(  \left(  \frac{3}{2}\right)
^{2/3}y_{k}-\frac{1}{k^{2/3}}+\frac{1}{9}\frac{\ln(k)}{k^{5/3}}\right)
=\frac{2C}{3},
\]%
\[
\lim_{k\rightarrow\infty}k^{1/3}\left(  \left(  \frac{3}{2}\right)
^{1/3}x_{k}-\frac{3}{2}k^{2/3}-\frac{1}{6}\frac{\ln(k)}{k^{1/3}}\right)
=C=0.8010888849039666437110775....
\]

\subsection{Case $q=4$}

We find%
\begin{align*}
4^{1/4}y_{k}  &  \sim\frac{1}{k^{1/4}}-\frac{3}{32}\frac{\ln(k)}{k^{5/4}%
}-\frac{C}{4}\frac{1}{k^{5/4}}+\frac{45}{2048}\frac{\ln(k)^{2}}{k^{9/4}%
}+\left(  -\frac{9}{256}+\frac{15C}{128}\right)  \frac{\ln(k)}{k^{9/4}}\\
&  +\left(  \frac{7}{256}-\frac{3C}{32}+\frac{5C^{2}}{32}\right)  \frac
{1}{k^{9/4}}-\frac{405}{65536}\frac{\ln(k)^{3}}{k^{13/4}}+\left(  \frac
{189}{8192}-\frac{405C}{8192}\right)  \frac{\ln(k)^{2}}{k^{13/4}}\\
&  +\left(  -\frac{297}{8192}+\frac{63C}{512}-\frac{135C^{2}}{1024}\right)
\frac{\ln(k)}{k^{13/4}}+\left(  \frac{75}{4096}-\frac{99C}{1024}+\frac
{21C^{2}}{128}-\frac{15C^{3}}{128}\right)  \frac{1}{k^{13/4}}%
\end{align*}
and thus
\[
-\lim_{k\rightarrow\infty}k^{5/4}\left(  4^{1/4}y_{k}-\frac{1}{k^{1/4}}%
+\frac{3}{32}\frac{\ln(k)}{k^{5/4}}\right)  =\frac{C}{4},
\]%
\[
\lim_{k\rightarrow\infty}k^{3/4}\left(  4^{3/4}x_{k}-4k^{1/4}-\frac{3}{8}%
\frac{\ln(k)}{k^{3/4}}\right)  =C=2.5097227971726886238611225....
\]

\subsection{Case $q=4/3$}

We find%
\begin{align*}
\left(  \frac{4}{3}\right)  ^{3/4}y_{k}  &  \sim\frac{1}{k^{3/4}}-\frac{3}%
{32}\frac{\ln(k)}{k^{7/4}}-\frac{3C}{4}\frac{1}{k^{7/4}}+\frac{21}{2048}%
\frac{\ln(k)^{2}}{k^{11/4}}+\left(  -\frac{3}{256}+\frac{21C}{128}\right)
\frac{\ln(k)}{k^{11/4}}\\
&  +\left(  \frac{5}{256}-\frac{3C}{32}+\frac{21C^{2}}{32}\right)  \frac
{1}{k^{11/4}}-\frac{77}{65536}\frac{\ln(k)^{3}}{k^{15/4}}+\left(  \frac
{27}{8192}-\frac{231C}{8192}\right)  \frac{\ln(k)^{2}}{k^{15/4}}\\
&  +\left(  -\frac{67}{8192}+\frac{27C}{512}-\frac{231C^{2}}{1024}\right)
\frac{\ln(k)}{k^{15/4}}+\left(  \frac{19}{4096}-\frac{67C}{1024}+\frac
{27C^{2}}{128}-\frac{77C^{3}}{128}\right)  \frac{1}{k^{15/4}}%
\end{align*}
and thus
\[
-\lim_{k\rightarrow\infty}k^{7/4}\left(  \left(  \frac{4}{3}\right)
^{3/4}y_{k}-\frac{1}{k^{3/4}}+\frac{3}{32}\frac{\ln(k)}{k^{7/4}}\right)
=\frac{3C}{4},
\]%
\[
\lim_{k\rightarrow\infty}k^{1/4}\left(  \left(  \frac{4}{3}\right)
^{1/4}x_{k}-\frac{4}{3}k^{3/4}-\frac{1}{8}\frac{\ln(k)}{k^{1/4}}\right)
=C=0.8248745112329031526004762....
\]

\subsection{Case $q=1/2$}

Starting with $z_{k+1}=z_{k}+\sqrt{z_{k}}$, $z_{0}=1$, let $x=\sqrt{z}$ and
obtain $x_{k+1}=\sqrt{x_{k}^{2}+x_{k}}$, i.e., the $p$-sequence with $p=1/2$.
\ The expansion of $x_{k}$ appears in Section 1.1; omitting details
(Proposition 12 of \cite{P1-popa} or Proposition 7 of \cite{P2-popa}), we
find
\begin{align*}
\lim_{k\rightarrow\infty}\frac{1}{k}\left(  z_{k}-\frac{1}{4}k^{2}+\frac{1}%
{4}k\ln(k)\right)   &  =1.1751774424585571398132856...\\
&  =0.5+0.6751774424585571398132856....
\end{align*}
Therefore we have come full circle. \ A\ name \textquotedblleft
add-the-square-root sequence\textquotedblright\ has been proposed for%
\[%
\begin{array}
[c]{ccccccccc}%
z_{0}=1, &  & z_{1}=2, &  & z_{2}=2+\sqrt{2}, &  & z_{3}=2+\sqrt{2}%
+\sqrt{2+\sqrt{2}}, &  & \ldots
\end{array}
\]
and the alternative constant appears merely because authors of \cite{SS-popa}
prescribed $\tilde{z}_{1}=1$, $\tilde{z}_{2}=2$, ... instead:\
\[
\tilde{z}_{k}=z_{k-1}\sim\frac{1}{4}(k-1)^{2}\sim\frac{1}{4}k^{2}-\frac{1}%
{2}k\sim z_{k}-\frac{1}{2}k
\]
hence $\tilde{z}_{k}/k\sim z_{k}/k-1/2$. \ Our search for a nontrivial
nonlinear recurrence with known $C$ (defined independently of the sequence,
e.g., an expression in closed-form) continues without success.

\section{Appendix}

The Mavecha-Laohakosol algorithm \cite{ML-popa} is applicable to an analytic
function $g(y)$ whose Taylor series at the origin starts as
\[
y+a_{1}y^{\tau+1}+a_{2}y^{2\tau+1}+a_{3}y^{3\tau+1}+\cdots
\]
where $a_{1}<0$ and $\tau\geq1$ is an integer. \ We demonstrate its use on the
reciprocal $p$-sequence $\{y_{k}\}$ with $p=1/2$, finding its asymptotic
expansion to order $1/k^{6}$. \ 

Formulaic knowledge of \cite{F4-popa} is assumed. \ Given $g(y)=y/\sqrt{1+y}$,
we have $\tau=1$,
\[
\{a_{m}\}_{m=1}^{7}=\left\{  -\frac{1}{2},\frac{3}{8},-\frac{5}{16},\frac
{35}{128},-\frac{63}{256},\frac{231}{1024},-\frac{429}{2048}\right\}
\]
and $\lambda=2$; consequently
\[
\{b_{j}\}_{j=1}^{6}=\left\{  -\frac{1}{2},\frac{1}{2},-\frac{5}{8},\frac{7}%
{8},-\frac{21}{16},\frac{33}{16}\right\}  ,
\]%
\[
\{a_{0j}\}_{j=1}^{6}=\left\{  1,-1,\frac{4}{3},-2,\frac{16}{5},-\frac{16}%
{3}\right\}  ,
\]%
\[
\{c_{i}\}_{i=1}^{5}=\left\{  0,\frac{1}{48},-\frac{1}{72},\frac{1}{240}%
,\frac{1}{225}\right\}  ,
\]%
\[
T_{2}=-\frac{1}{2}Y,
\]%
\[
T_{3}=\dfrac{1}{4}Y^{2}+\dfrac{1}{4}Y-\dfrac{1}{48},
\]%
\[
T_{4}=-\dfrac{1}{6}Y^{3}-\frac{3}{8}Y^{2}-\dfrac{1}{12}Y+\dfrac{7}{288},
\]

\[
T_{5}=\frac{1}{8}Y^{4}+\frac{11}{24}Y^{3}+\frac{5}{16}Y^{2}-\frac{1}%
{32}Y-\frac{47}{2880},
\]%
\[
T_{6}=-\frac{1}{10}Y^{5}-\frac{25}{48}Y^{4}-\frac{31}{48}Y^{3}-\frac{3}%
{32}Y^{2}+\frac{233}{2880}Y+\frac{41}{14400}%
\]
and
\[
P_{2}=Y^{2}-\frac{1}{2}Y,
\]%
\[
P_{3}=Y^{3}-\dfrac{5}{4}Y^{2}+\dfrac{1}{4}Y+\dfrac{1}{48},
\]%
\[
P_{4}=Y^{4}-\dfrac{13}{6}Y^{3}+\dfrac{9}{8}Y^{2}-\dfrac{1}{24}Y-\dfrac{7}%
{288},
\]%
\[
P_{5}=Y^{5}-\frac{77}{24}Y^{4}+\dfrac{71}{24}Y^{3}-\dfrac{2}{3}Y^{2}%
-\dfrac{29}{288}Y+\dfrac{47}{2880},
\]%
\[
P_{6}=Y^{6}-\frac{87}{20}Y^{5}+\frac{145}{24}Y^{4}-\dfrac{45}{16}Y^{3}%
+\dfrac{1}{32}Y^{2}+\dfrac{427}{2880}Y-\dfrac{139}{57600}.
\]
The remarkable formula connecting $P_{m}$ and asymptotics of $y_{k}%
=g(y_{k-1})$ is%
\[
y_{k}\sim\left(  \frac{\lambda}{k}\right)  ^{1/\tau}\left\{  1+%
{\displaystyle\sum\limits_{m=1}^{6}}
P_{m}\left(  -\frac{1}{\tau}\left[  b_{1}\ln(k)+2C\right]  \right)  \frac
{1}{k^{m}}\right\}  ,
\]
the first four terms of which imply%
\begin{align*}
y_{k}  &  \sim\frac{2}{k}+\frac{\ln(k)}{k^{2}}-\frac{4C}{k^{2}}+\frac{1}%
{2}\frac{\ln(k)^{2}}{k^{3}}-\left(  \frac{1}{2}+4C\right)  \frac{\ln(k)}%
{k^{3}}+\left(  2C+8C^{2}\right)  \frac{1}{k^{3}}+\frac{1}{4}\frac{\ln(k)^{3}%
}{k^{4}}\\
&  -\left(  \frac{5}{8}+3C\right)  \frac{\ln(k)^{2}}{k^{4}}+\left(  \frac
{1}{4}+5C+12C^{2}\right)  \frac{\ln(k)}{k^{4}}-\left(  -\frac{1}{24}%
+C+10C^{2}+16C^{3}\right)  \frac{1}{k^{4}}%
\end{align*}
(echoing Section 1.1). \ The fifth \&\ sixth terms substantially extend the
series:
\begin{align*}
&  \frac{1}{8}\frac{\ln(k)^{4}}{k^{5}}-\left(  \frac{13}{24}+2C\right)
\frac{\ln(k)^{3}}{k^{5}}+\left(  \frac{9}{16}+\frac{13}{2}C+12C^{2}\right)
\frac{\ln(k)^{2}}{k^{5}}\\
&  -\left(  \frac{1}{24}+\frac{9}{2}C+26C^{2}+32C^{3}\right)  \frac{\ln
(k)}{k^{5}}+\left(  -\frac{7}{144}+\frac{1}{6}C+9C^{2}+\frac{104}{3}%
C^{3}+32C^{4}\right)  \frac{1}{k^{5}}\\
&  +\frac{1}{16}\frac{\ln(k)^{5}}{k^{6}}-\left(  \frac{77}{192}+\frac{5}%
{4}C\right)  \frac{\ln(k)^{4}}{k^{6}}+\left(  \frac{71}{96}+\frac{77}%
{12}C+10C^{2}\right)  \frac{\ln(k)^{3}}{k^{6}}\\
&  -\left(  \frac{1}{3}+\frac{71}{8}C+\frac{77}{2}C^{2}+40C^{3}\right)
\frac{\ln(k)^{2}}{k^{6}}+\left(  -\frac{29}{288}+\frac{8}{3}C+\frac{71}%
{2}C^{2}+\frac{308}{3}C^{3}+80C^{4}\right)  \frac{\ln(k)}{k^{6}}\\
&  -\left(  -\frac{47}{1440}-\frac{29}{72}C+\frac{16}{3}C^{2}+\frac{142}%
{3}C^{3}+\frac{308}{3}C^{4}+64C^{5}\right)  \frac{1}{k^{6}}.
\end{align*}

Our simple procedure for estimating the constant $C$ involves computing
$y_{K}$ exactly via recursion, for some suitably large index $K$. \ We then
set the value $y_{K}$ equal to our series and numerically solve for $C$. \ The
effect of using the additional series terms is fairly dramatic. \ When
employing terms up to order $1/k^{4}$, an index $\approx10^{10}$ may be
required for $25$ digits of precision in the $C$ estimate. \ By way of
contrast, when employing terms up to order $1/k^{6}$, an index $\approx10^{6}$
might suffice.

\section{Addendum}

The $q$-sequence with $q=2$, when squared, has asymptotics%
\[
x_{k}^{2}\sim2k+\frac{1}{2}\ln(k)+2C.
\]
If we prescribe $\tilde{x}_{1}=1$, $\tilde{x}_{2}=2$, ... instead:%
\[
\tilde{x}_{k}^{2}=x_{k-1}^{2}\sim2(k-1)\sim x_{k}^{2}-2
\]
then%
\[
\tilde{x}_{k}^{2}\sim2k+\frac{1}{2}\ln(k)+2(C-1)
\]
and our estimate%
\[
2(C-1)=-0.2768576248625765389364372...
\]
improves upon a figure $-0.2768576$ given years ago by Jean Anglesio
\cite{BL-popa}.

The same sequence $\{x_{k}\}$ arises in graph theory \cite{LPD-popa},
arithmetic dynamics \cite{Slv-popa} and approximation \cite{Mnt-popa,
F6-popa}. \ Define an infinite product%
\[
r=%
{\displaystyle\prod\limits_{k=0}^{\infty}}
\left(  1+\frac{1}{x_{k}^{2}}\right)  ^{1/2^{k+1}}=1.54170091336287603176....
\]
Consider the class $S_{n}$~of all integer polynomials of the exact degree $n$
and all $n$ zeroes both in $[-1,1]$ and simple.\ Let%
\[%
\begin{array}
[c]{ccccc}%
{\displaystyle\sum\limits_{k=0}^{n}}
a_{k,n}x^{k}\in S_{n}, &  & a_{n,n}\neq0, &  & n=1,2,3,\ldots
\end{array}
\]
be an arbitrary sequence $R$ of polynomials. \ It is known that
\cite{Pri-popa}%
\[
\inf_{R}\operatorname*{liminf}\limits_{n\rightarrow\infty}\left\vert
a_{n,n}\right\vert ^{1/n}\leq r
\]
and quite possible that no better bound exists. \ It is fascinating that
$x_{k}+1/x_{k}$ should have occurred \textit{at all} in these problems. \ We
speculate about other unexpected places where nonlinear recurrences might emerge.

Just as $q=1/2\longleftrightarrow p=1/2$, it can be shown that
$q=2\longleftrightarrow p=2$. \ Starting with $z_{k+1}=z_{k}+1/z_{k}$,
$z_{0}=1$, let $x=z^{2}$ and obtain%
\[
\sqrt{x_{k+1}}=\sqrt{x_{k}}+\frac{1}{\sqrt{x_{k}}}%
\]
i.e.,
\[
x_{k+1}=x_{k}+2+\frac{1}{x_{k}}=\frac{x_{k}^{2}+2x_{k}+1}{x_{k}}%
=\frac{(1+x_{k})^{2}}{x_{k}}%
\]
i.e., the $p$-sequence with $p=2$. \ Note that
$1.723142375137...=(2)(0.861571187568...)$ as a consequence.

\section{Acknowledgements}

I am grateful to Popa for a very helpful discussion. \ In the statement of
Theorem 5 of \cite{P2-popa}, the lead coefficient $1/2$ of the $1/n$ term
should be $1/4$. \ Also, in the statement of Theorem 10 of \cite{P2-popa}, the
lead coefficient $2a_{1}^{4}K^{2}$ of the $1/n^{3}$ term should be $2a_{1}%
^{6}K^{2}$ (Popa's corresponding expression $2a_{1}^{2}C^{2}$ on the following
page, where $C=-a_{1}^{2}K$, is correct). \ This algebraic error affects his
Corollary 11. \ The creators of Mathematica earn my gratitude every day:\ this
paper could not have otherwise been written.

\end{document}